\documentclass[12pt]{article} 
\usepackage[sectionbib]{natbib}
\usepackage{array,epsfig,fancyheadings,rotating,todonotes}
\usepackage[]{hyperref}  
\usepackage{sectsty, secdot}
\sectionfont{\fontsize{12}{14pt plus.8pt minus .6pt}\selectfont}
\subsectionfont{\fontsize{12}{14pt plus.8pt minus .6pt}\selectfont}
\usepackage{multirow}
\usepackage{soul}
\usepackage{amsmath}
\usepackage{amssymb}
\usepackage{amsfonts}
\usepackage{multirow}
\usepackage{amsthm}
\usepackage{booktabs}
\usepackage{cleveref}

\newtheorem{theorem}{Theorem}[section]
\newtheorem{lemma}[theorem]{Lemma}
\newtheorem{remark}{Remark}
\newtheorem{proposition}{Proposition}
\newtheorem{corollary}{Corollary}
\theoremstyle{remark}
\newtheorem{definition}[theorem]{Definition}

\usepackage{amssymb,amsmath}

\usepackage{graphicx,color}
\usepackage[utf8x]{inputenc}
\usepackage{amsmath}
\usepackage{amsfonts}
\usepackage{amssymb}
\usepackage{amsthm}
\usepackage{graphicx}
\usepackage{array}
\usepackage{multirow}
\usepackage{fancyhdr}
\usepackage{setspace}
\usepackage{makeidx}
\usepackage{hyperref}

\newcommand{\eps}{\varepsilon}

\newcommand{\ind}{\mathbb{I}}

\newcommand{\prob}{\mathbb{P}}
\newcommand{\normal}{\mathbf{n}}

\newcommand{\gen}{\mathcal{L}}

\begin{document}
	
	
	\begin{center}
		\large \bf   Home range estimation under a restricted sampling scheme.\normalsize
	\end{center}
	\normalsize

	\begin{center}
		Alejandro Cholaquidis$^a$, Ricardo Fraiman$^a$,  and Manuel Hernandez-Banadik$^b$\\
		
		$^a$ Centro de Matem\'atica, Universidad de la Rep\'ublica, Uruguay\\
		$^b$ Instituto de Estadística, Universidad de la Rep\'ublica, Uruguay 
		\\
	\end{center}
	\begin{abstract}
		
		The analysis of animal movement has gained attention recently. New continuous-time models and statistical methods have been developed to estimate some sets related to their movements, such as the home-range and the core-area among others, when  the information of the trajectory is provided by a GPS.
		Because data transfer costs and GPS battery life are practical constraints in ecological studies, the experimental designer must make critical sampling decisions in order to maximize information. To capture fine-scale motion, long-term behavior must be sacrificed, and vice versa.
		To overcome this limitation, we introduce the on--off sampling scheme,  where the GPS is alternately on and off. This scheme is already used in practice but with insufficient statistical theoretical support.
		We prove the consistency of home-range estimators with an underlying reflected diffusion model under this sampling method (in terms of the Hausdorff distance). The same rate of convergence is achieved as in the case where the GPS is always on for the whole experiment. This is illustrated by a simulation study and real data. We also provide estimators of the stationary distribution, its level sets (which give estimators of the core-area), and the drift function.
		
	\end{abstract}

	\section{Introduction}
	
	Home-range was  defined by \cite{burt:43} as ``the area traversed by the individual in its normal activities of food gathering, mating, and caring for young.'' 
	Its estimation is a major problem in animal ecology.
	Hence there has arisen a considerable literature on the subject, see for instance the reviews in \cite{worton:87}, \cite{powell:00} or \cite{chacon21}.
	Another key concept is the so-called core-area (areas where animals spend most of their time); in our setup, core-areas can be modelled by the level sets of the stationary distribution.
	Several models have been proposed to analyse the home-range as well as the dynamics of animals.
	The first  one \citep{hayne:49} assumed that the available data are a sequence of locations recorded at some times, and estimates the home-range by means of the convex hull of the points.
	More flexible  proposals have also  been introduced, for instance the use of ``local convex hulls'' \citep{getz:04} or the $r$-convex hull \citep{burgman:03}, since $r$-convex sets are a much more flexible class of sets.
	Some others,  referred to as occurrence estimators,  are focused on the estimation of the so-called ``utilization distribution'' (the density function that describes the probability of finding the animal at a particular location).
	Among these, there is the Brownian bridge model (BBMM),  a parametric model introduced in \cite{horne}.

	The importance of using continuous time models for animal movement is discussed in  \cite{kie, no, fleming2014fine, fleming:15, calabrese2016ctmm}.
	Home-range estimation is carried out by estimating the support of a probability density function, or an appropriate level set.
	The  importance of taking into account the autocorrelations is deeply analysed in those papers.
	They reproduce different autocorrelation scenarios (position or velocity autocorrelation) by using continuous-time movement models and contrast a kernel-based home-range estimator designed for iid data, with the autocorrelated kernel density estimation (AKDE), tailored for autocorrelated data \citep{fleming:15}.
	Their findings confirm that the former  underestimates the home-range, while the latter one performs very well.
	
	AKDE adjusts the smoothing kernel parameter in order to weight properly the contribution of each location point to the estimation of the density: the more dependent the data are, the less `new information' they provide.
	So, AKDE aims to fit an autocovariance function to compute the estimator.
	To this end, several animal movement models are proposed, such as the Ornstein--Uhlenbeck, the Ornstein--Uhlenbeck with foraging \citep[Supplementary matherial C]{fleming2014fine}, and integrated process are also considered (that ensure velocity autocorrelation), among others (see \cite{no, fleming2014fine, fleming:15, calabrese2016ctmm}).
	From a statistical point of view, no asymptotic theory has yet been developed and, more importantly, some of the aforementioned models lack a stationary distribution to be estimated.

	From a more theoretical-mathematical point of view, in \cite{ch:16} the animal's motion is modelled as a reflected Brownian motion on a bounded set $S \subset \mathbb{R}^d$ that plays the role of the home-range, and consistent estimators and rates of convergence are obtained.
	In \cite{ch:20} a reflected Brownian motion with drift is considered (RBMD, see section \ref{crof}), and it is proven that, under suitable conditions, there exists a unique stationary distribution of the process, and consistent uniform estimators are obtained for its density, together with its level sets (that give an estimate of the core-area).

 Even though it seems unrealistic to model animal movement with an RBMD (due to the non-regularity of the paths), since the data are stored in a computer as a vector of locations,  it can not be distinguished if the underlying process is smooth or not.
	The choice of the RMBD is  mainly due to its good properties (it is a geometrically ergodic Markov process), and it comes as a theoretical tool. 
	Indeed, all the results we present remain true if we assume that the discrete vector of locations is a geometrically ergodic Markov process in discrete time, while as suggested in \cite{no} our model takes into account the position autocorrelation.
	The RBMD provides an example of a stochastic process that fulfills these conditions, and provides theoretical support for some previous empirical results, such as those in \cite{no}.

	A problem that appears in practice is that these models assume that it is possible to record the location over a long period of time to attain convergence, but battery life is an actual constraint that must be taken into account.
	A hands-on way researchers found to overcome this limitation was decreasing the sampling frequency of GPS device, so the duration of the sampling gets better.
	The trade-off between battery life and fine-scale movement resolution is of paramount importance, as mentioned in \cite{brown12},
	``if intervals are too long, they undersample the details of movement paths, and if too short, they oversample resting sites and deplete the unit’s battery without providing new information."

	\cite{brown12, brown13} propose combining GPS information with information from an accelerometer, which automatically turns the GPS off when the animal is not moving.
	
	So, for this purpose, several estimators have been proposed and evaluated both in empirical and simulation studies, but there has not been much about how the data were collected except the sampling frequency because of the autocorrelation introduced. As is pointed in \cite{mitchell2019trade}: ``Research on home-range and habitat selection for these species should therefore incorporate a consideration of how different sampling parameters and methods may affect the structure of the data and the conclusions drawn.
	However, factors such as these are seldom explicitly considered.''
	
	This issue is recently faced in \cite{he2022guide}: ``using GPS devices requires making a number of decisions about sampling that can affect the robustness of a study's conclusion''; they provide an exhaustive practitioner-oriented anaylsis for designing GPS-based tracking studies.
	It focuses on partitioning the sampling effort over time: ``it can be useful to vary GPS sampling across time to focus data collection on certain time windows over others.
	[...] GPS devices are often switched off at night to conserve battery.''
	Several practical factors must be considered: because of the weight of the battery, small animals cannot be sampled for long time periods as bigger animals can be; furthermore, ``battery consumption does not scale linearly with increasing sampling frequency.
	GPS devices use energy each time they (re)boot and search for satellites, which typically takes longer with sparser sampling regime'' \citep{he2022guide}; there are solar-powered GPS devices that in some cases take a fixed number of days to recharge their batteries, so the sampling schedule becomes periodic (see \cite{he2022guide} and references therein).

\begin{figure} 
		\begin{center}
			\includegraphics[width=.85\textwidth]{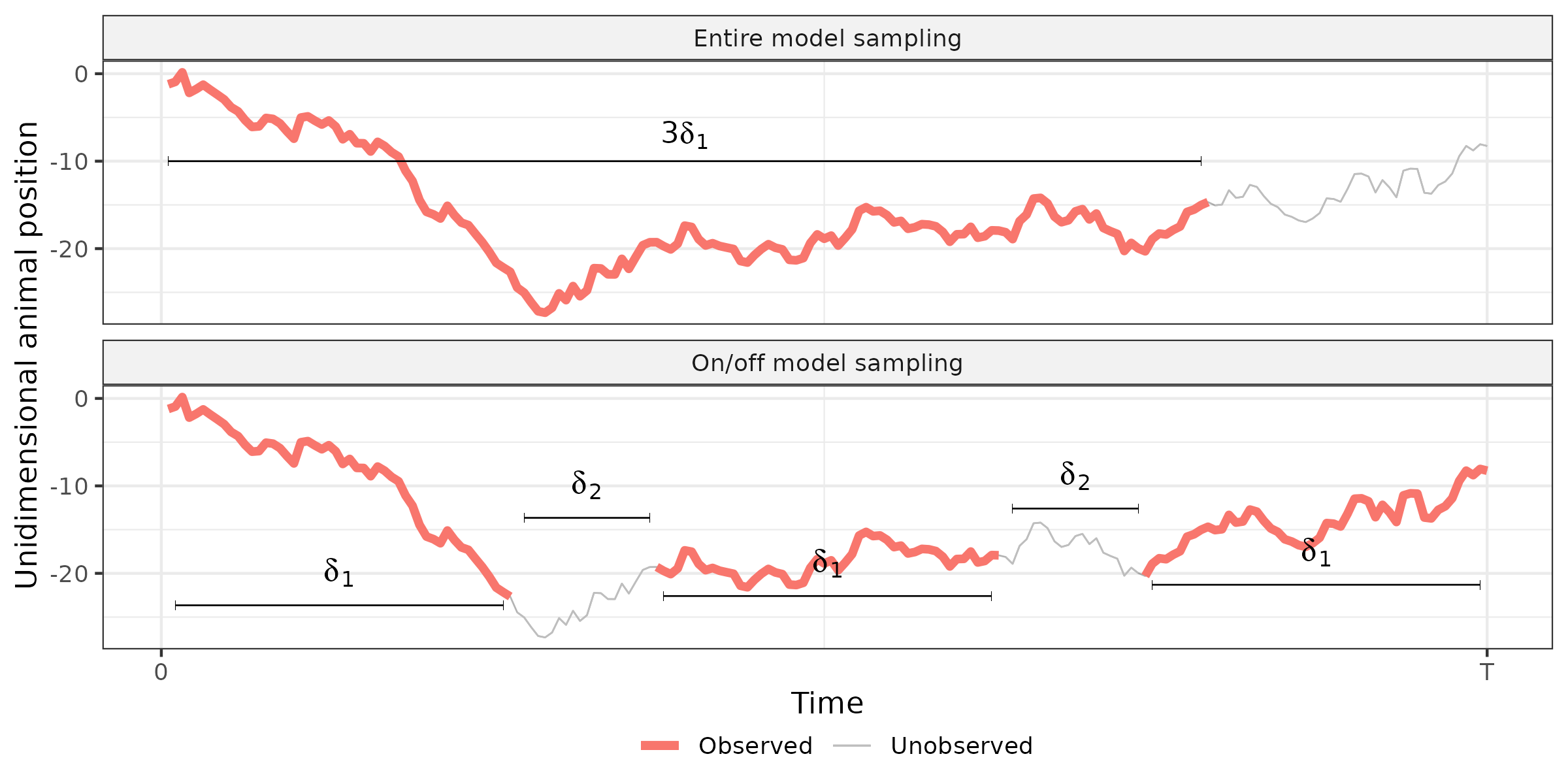} 
			\caption{Top panel: in red a one dimensional process observed from time 0 to time $3\delta1$.
				Bottom panel: the same process but observed intermittently, the trajectory where the GPS is on is shown in red while in grey is shown the unobserved trajectory where the GPS is off.}
			\label{fig:modelo_observacion}
		\end{center}
	\end{figure}

	So, more complex sampling schemes have been introduced in the literature but only, as far as our knowledge extends, from a practical perspective and no theoretical results support their findings.
	We propose a sampling method that consists of recording the position of the animal continuously, but only for certain periods of time in which the GPS is on (see Figure \ref{fig:modelo_observacion} for an example of applying this sample scheme to a simulated unidimensional trajectory).
	More precisely, the GPS is set alternately on during $p$ intervals of length $\delta_1$, and off during $p-1$ intervals of length $\delta_2$.
	This will be called the on--off model (see Section \ref{sec-nutshell}).
	The existence and uniqueness of such a stochastic process obtained by means of this sampling scheme is inherited from the RBMD itself.
	It might be considered as an intermediate solution of the trade-off explained above: it allows a higher sampling frequency but also improves the use of the battery life.
	
	The intuition behind this is that statistical properties should not vary so much when observing the full trajectory compared to when it is observed intermittently.
	Our findings confirm  the empirical rule (see Eq. (5) in \cite{no}) which states that ``the information content of a tracking data set is not a function of the total number of locations, [...] but rather the equivalent number of statistically independent locations'', which is usually called the effective sample size.

	Assuming the RBMD underlying model,  and denoting by $T$ the whole time where the process is observed (i.e $T=p\delta_1+(p-1)\delta_2$), we show that any set $S_T^{\textrm{ON}}\subset S$ containing the trajectory of the on--off model  is a consistent estimator of $S$ (in the Hausdorff distance) when  $T\to\infty$
	We show that the rate of convergence of the estimator is the same as when the full trajectory is available, as $T\to\infty$ and $\delta_1+\delta_2\to \infty$ appropriately.
	The obtained rate, for $S \subset \mathbb R^2$, $\mathcal{O}(\log(T)^2/T)$, is close to the rate when the data is an iid sample of $n$ points in $S$ (which is proven to be sharp in \cite{cuevas1997plug}).
	As often  happens in several nonparametric estimation problems, when we change from an iid setting to a mixing one, an extra $\log$ term appears.

	Next we consider the case where $S$ is $r$-convex (see the definition in the next section) and provide  consistency results  with respect to the Hausdorff distance and to the distance in measure for the $r$-convex hull.
	The notion of $r$-convex sets has been extensively studied in various areas, particularly in stochastic geometry and set estimation, since it provides a much more general and flexible family of sets than that of convex sets, maintaining some important properties.
	We also provide estimators of the stationary distribution, level sets and the drift function.
	Lastly we provide some simulation results and a real data example is analysed.

	\section{Some Preliminaries}
	
	\noindent\textit{The following notation will be used throughout this paper}.\\
	
	Given a set $S\subset \mathbb{R}^d$, we will denote by
	$\textnormal{int}(S)$, $\overline{S}$,  $\partial S$ and \textcolor{red}{$S^c$} the interior, closure, boundary and complement of $S$,
	respectively, with respect to the usual topology of $\mathbb{R}^d$.

	The parallel set of $S$ of radius $\varepsilon$ will be denoted by $B(S,\eps)$, that is,
	$B(S,\eps)=\{y\in{\mathbb R}^d:\ \inf_{x\in S}\Vert y-x\Vert\leq \eps \}$.
	If $A\subset\mathbb{R}^d$ is a Borel set, then $\mu(A)$ denotes
	its $d$-dimensional Lebesgue measure.
	We will denote by $B(x,\eps)$ the closed ball
	in $\mathbb{R}^d$,
	of radius $\eps$, centred at $x$, and  $\omega_d=\mu(B(0,1))$.
	The open ball is denoted by $\mathring{B}(x,\eps)$.
	Given two compact non-empty sets $A,C\subset{\mathbb R}^d$, 
	the Hausdorff distance  between $A$ and $C$ is defined by
	\begin{equation*}
		d_H(A,C)=\inf\{\eps>0: \mbox{such that } A\subset B(C,\eps)\, \mbox{ and }
		C\subset B(A,\eps)\}.
	\end{equation*}
	Given two measurable sets $A, C \subset \mathbb{R}^d$, the distance in measure between $A$ and $C$ is defined by 
	\begin{equation*}
		d_\mu(A, C) = \mu(A\setminus C)+\mu(C\setminus A).
	\end{equation*}
	Next, we introduce  the  $r$-convex sets, a well-known shape restriction in set estimation  \citep{walther:97,walther:99},  which extends  convex sets to a much more flexible family of sets.
	It just replaces the hyperplanes in the definition of convex sets by the complements of balls of radius $r$, providing a very flexible class of sets.
	
\begin{figure}[h]
		
		\includegraphics[width=.85\textwidth]{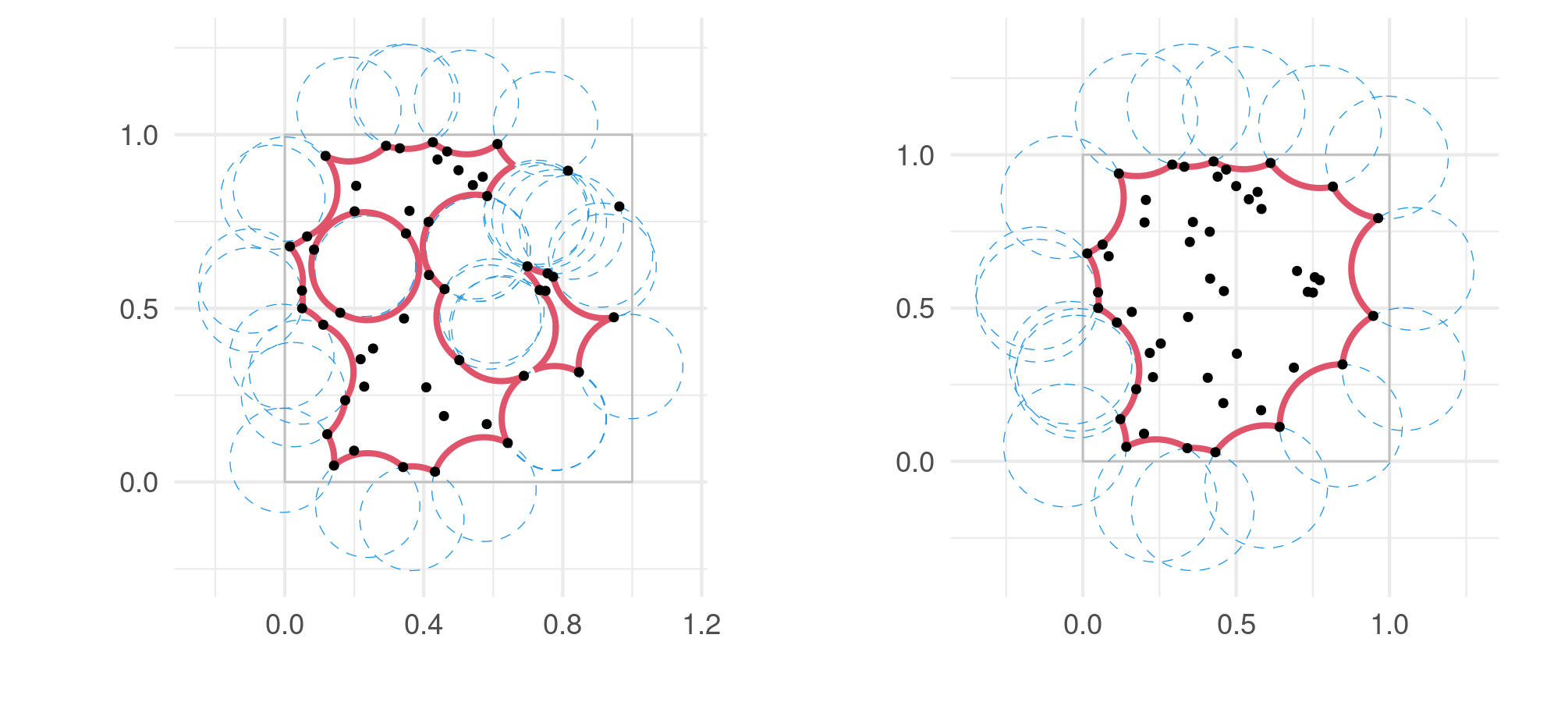}
		\caption{In solid red line the $r$-convex hull of a sample of 100 uniformly distributed random vectors in $[0,1]^2$.
			At left $r = 0.15$, while at right $r = 0.2$.}\label{fig:cierre-r-convexo}
		
	\end{figure}

	\begin{definition} A set $S\subset \mathbb{R}^d$ is said to be $r$-convex, for $r>0$, if 
		$S=C_r(S),$ where
		\begin{equation*} 
			C_r(S)=\bigcap_{\big\{ \mathring{B}(x,r):\ \mathring{B}(x,r)\cap S=\emptyset\big\}} \Big(\mathring{B}(x,r)\Big)^c,
		\end{equation*}
		is the $r$-convex hull of $S$.	
	\end{definition}
	Refer to Figure \ref{fig:cierre-r-convexo} for an example of the $r$-convex hull of a set of points for two different values of  $r$.
	
	Let us recall the definition of a $C^2$ boundary. Intuitively, this means that the boundary, locally, can be thought of as the graph of a $C^2$ function.
	
	\begin{definition}[From \cite{evans2010partial}]
		Given a bounded set $S \subset \mathbb{R}^d$, we say that $\partial S$ is $C^2$ if for every point $x \in \partial S$, there exists  $r >0$ and a function of class $C^2$, $\gamma: \mathbb{R}^{d-1} \to \mathbb{R}$, such that --upon relabelling and reorienting axes-- we have $$S \cap B(x, r) = \{(z_1,\dots, z_d) \in B(x, r): z_d > \gamma(z_1, \dots, z_{d-1})\}.$$	
	\end{definition}

	\subsection{A brief outline of our theoretical results}
	
	We  prove that any set $S_T^{ON} \subset B(S,\eps_T)$ for $\eps_T \to 0$ containing the trajectory of the on--off model (that is, the locations at which the GPS is on) is a consistent estimator, in the Hausdorff distance, of the home-range.
	
	However, this is not the case if we want to estimate the set with respect to the distance in measure (w.r.t. Lebesgue measure).
	In this case, we propose using  the $r$-convex hull of the same trajectory of the on--off model.
	When $S$ is $r$-convex, a natural estimator of $S$ from a random sample $\aleph_n$ of points 
	(drawn from a distribution with support $S$) is $C_r(\aleph_n)$, \citep{rodriguez:07,pateiro09}.
	We show its convergence in the Hausdorff metric, while the convergence in measure is derived from Corollary \ref{coro}, as mentioned in Remark \ref{rem2}.
	
	In order to estimate the core-area   and the drift component of the stochastic differential equation (\ref{sde}), we prove the uniform convergence of a kernel estimator of the stationary density, and derive estimators of the level sets in Theorem \ref{teo:levels}.
	
	An estimator of the the drift component is derived from the kernel density estimator. 
	
	Some techniques used in the proofs resemble the ones used in  \citep{ch:16,ch:20}.
	However, the rate of convergence of the Hausdorff distance for the RBMD is not provided in the aforementioned works, nor in the literature, as far as  our knowledge extends, and the behaviour of the on--off model is not considered.
	
	\subsection{Reflected Brownian motion with drift}\label{crof}
	
	Now we will give a brief review of the definition and main properties of the RBMD.
	The details can be found
	for instance in \cite{ch:20} and references therein.
	In what follows, ${S}\subset \mathbb{R}^d$ is such that $S=\overline{\textnormal{int}(S)}$ and $\textnormal{int}(S)$ is a bounded domain in $\mathbb{R}^d$ (that is, a bounded connected open set) such that $\partial S$ is $C^2$.
	Given a $d$-dimensional Brownian motion $\{B_t\}_{t\geq 0}$ departing from 
	$B_0=0$  defined on a filtered probability space $(\Omega,\mathcal{F},\{\mathcal{F}_t\}_{t\geq 0},\prob)$, 
	and a function $f:S\rightarrow \mathbb{R}$,  the RBMD is the (unique) strong solution to the following reflected stochastic differential equation on $S$ whose drift is given by 
	the gradient of  $f$,
	\begin{equation}\label{sde}
		X_t= X_0+ B_t-\frac{1}{2}\int_0^t\nabla f(X_s)ds+\int_0^t\normal(X_s)d\xi_s,
		\quad\text{ where } X_t\in S,\ \forall t\geq 0.
	\end{equation}
	Here we assume that $\nabla f$ is Lipschitz, 
	while $\normal(x)$ denotes the inner unit vector at the boundary point $x\in\partial S$.
	The term $\{\xi_t\}_{t\geq 0}$ is the corresponding \emph{local time}, that is, a one-dimensional continuous non-decreasing process with $\xi_0=0$ that satisfies
	\[
	\xi_t=\int_0^t\mathbb{I}_{\{X_s\in\partial S\}}d\xi_s,
	\]
	see \cite{saisho:87}.
	Since we have assumed that $\partial S$ is $C^2$,  we can ensure that the geometric conditions for the existence of a solution of Equation \eqref{sde}, as required in \cite{saisho:87}, are satisfied.
	We then get from Theorem 5.1 in \cite{saisho:87} that there exists an unique strong solution of the Skorokhod stochastic differential equation \eqref{sde}.
	The solution is a strong solution in the sense of Definition 1.6 in \cite{ikeda}.
In \cite{ch:20} it is proved that such a set $S$ as mentioned above has an invariant distribution.  Recall that a probability measure $\pi$ on $S$ is said to be an invariant measure for a time-homogeneous Markov process $\{Z_t\}_{t\geq 0}$ if $\int_S \mathbb{P}_x(Z_t\in A)\pi(dx)=\pi(A)$ for all $t>0$ and all Borel set $A\subset S$, where $\mathbb{P}_x$ denotes the probability assuming $X_0=x$ a.s.

	The following proposition, proven in \cite{ch:20}, will be used to get the consistency in Hausdorff distance of the trajectory, as an estimator of the home-range of the on--off model.
	
	\begin{proposition} \label{properg} 
		Let $S\subset \mathbb R^d$ be such that $\textnormal{int}(S)$ is  a bounded domain and  $\partial S$ is $C^2$.
		Denote by $\pi$ the invariant distribution of $\{X_t\}_{t\geq 0}$. Then 
		there exist positive constants $\alpha$ and $\beta$ such that 
		\begin{equation*}
			\sup_{x \in \textnormal{int}(S)} \big\| \mathbb P_x (X_t \in \cdot ) -\pi(\cdot)\big\|_{TV} \leq \beta e^{-\alpha t}.
		\end{equation*}
		Here, $\|\cdot\|_{TV}$ stands for the total variation norm of a measure, and $\mathbb{P}_x$ denotes the probability assuming $X_0=x$ a.s.
	\end{proposition}
	
	\begin{remark} The process defined by \eqref{sde} can be seen as a limit of a random walk sampling, that is, at each time we sample a point and a direction and magnitude to move.
		A precise proof of this is given in \cite{burdzy2008discrete} and \cite{bossy2004symmetrized}.
		This  sampling scheme is used in our simulations.
	\end{remark}

	\section{The On--Off Model}\label{sec-nutshell}
	Our on--off model is defined as follows:
	\begin{definition} 
		Given 
		\begin{itemize}
			\item $S\subset \mathbb{R}^d$ a compact set.
			\item $\{X_t:t>0\}$ a reflected Brownian motion with drift, in $S$.
			\item Two parameters $\delta_1, \delta_2 \in \mathbb{R}^+$.
			\item A function $\{a_t:t>0\}$  that varies over $\{0,1\}$ intermittently, with $a_t=1$ for periods of length $\delta_1$ and $a_t=0$ for periods of length $\delta_2$.\\
			More precisely, 
			$$a_t = \sum_{k=0}^\infty \ind_{[k(\delta_1 + \delta_2), (k+1)\delta_1 + k\delta_2]}(t).$$
		\end{itemize} We define the process  
		\begin{equation} \label{modelo}
			X_T^{\mathrm{ON}}=\{X_t : t \in \mathcal{I}, t<T \},
		\end{equation} 
		where $\mathcal{I}:=\{t: a_t=1\}$.
		Observe that the process $X_T^{\mathrm{ON}}$ is defined only on a  union of disjoint intervals.
		The function $a_t$ works like an on--off switch: we  only observe the process while $a_t=1$ (i.e. the switch is `On'), which happens on intervals of length $\delta_1$, while it is not observable  on intervals of length $\delta_2$, in alternation.
	\end{definition}

	The following theorem provides the almost sure convergence of the Hausdorff distance between any set, $S_T^{\mathrm{ON}}\subset S$ containing the trajectory of the on--off model and the set $S$, when $T\to \infty$ for fixed $\delta_1$ and $\delta_2$.
	If   $\delta_1+\delta_2\to \infty$, the same rate of convergence is obtained as in the case where the process $X_t$ is always observed over the whole experiment.
	Also, the obtained rate are close to $(n/\log(n))^{1/d}$, which was shown to be sharp in the iid case, that is, when instead of a process we have an iid sample of $n$ points in $S$ (see \cite{cuevas1997plug}).

	\begin{theorem}  \label{onoffTEO} 	Let $S \subset \mathbb{R}^d$ be a compact set such that $S = \overline{\textnormal{int}(S)}$ and $\partial S$ is $C^2$.
		Let $X_T^{\mathrm{ON}}$ be defined as in Equation (\ref{modelo}).
		Suppose that the drift is a Lipschitz function given by the gradient of some function $f$, assume that the stationary distribution $\pi$ has density $g$, 
		and that $c:=\inf_{x\in S} g(x)>0$.

		\begin{itemize}
			\item [a)] If $\delta_1,\delta_2>0$ are fixed,
			\begin{equation*}
				d_H(S_T^{\mathrm{ON}},S)\to 0\quad \mathrm{a.s.},
			\end{equation*}   
			for any set $S_T^{\mathrm{ON}}\subset S$ containing $X_T^{\mathrm{ON}}$ a.s.
			\item[b)] Let $\kappa_T= (T/\log^2(T))^{1/d}$. If 
			$$\delta_1 +\delta_2=-\frac{1}{\alpha}\log\Bigg(\frac{c}{2\beta }\Bigg(\frac{1}{\kappa_T\eta_T}\Bigg)^d\Bigg),$$
			where $\eta_T$ is any sequence that tends to infinity, then
			\begin{equation*}
				\kappa_Td_H(S_T^{\mathrm{ON}},S)\to 0\quad \mathrm{a.s.}
			\end{equation*}   
			for any set $S_T^{\mathrm{ON}}\subset S$ containing $X_T^{\mathrm{ON}}$ a.s.
		\end{itemize}
	\end{theorem}

	\begin{remark} \label{losdelta} The choice of the parameters is an important practical problem.
		In practice, the battery life, $B=p\delta_1$, and the time $T$ of the experiment remain fixed.
		Studying the optimization problem
		$$\min_{\{\delta_1,\delta_2:T=B/\delta_1+(B/\delta_1-1)\delta_2\}} d_H(S_T^{\mathrm{ON}},S)$$
		for $T$ and $B$ fixed,	is far beyond the scope of this paper and remains as an open problem.
		However,  we consider different alternatives in our simulations that confirm the empirical rule known in the literature of home-range estimation, which states that decreasing the correlation in the dataset increases the effective sample size \citep{no}.  Regarding the choice of $S_T^{\mathrm{ON}}$, if, for instance, $S$ is assumed to be $r$-convex or $\rho$-cone-convex (see Definition 2 in \cite{ch:14}), the $r$-convex hull (or the $\rho$-cone convex hull) of $X_T^{\textrm{ON}}$ can be used  as $S_T^{\textrm{ON}}$.
	\end{remark}

	The following corollary is a direct consequence of Theorem \ref{onoffTEO}, using  that if $\partial S$ is $C^2$, then it is $r$-convex for some $r>0$, and then  $C_r(X_T^{\mathrm{ON}})\subset S$ a.s.
	
	\begin{corollary} \label{coro} Under the hypotheses of Theorem \ref{onoffTEO}, for any measurable set  $S_T^\mathrm{ON}\subset S$ containing  $X_T^{\mathrm{ON}}$ {a.s.}, we have  that, for some $r>0$,
		$\kappa_Td_H(C_{r}(S_T^{\mathrm{ON}}),S)\to 0$ \ a.s., as $T\to \infty$, when $\delta_1,\delta_2$ are as in Theorem \ref{onoffTEO}.
	\end{corollary}
	
	\begin{remark} \label{rem2} The convergence in Hausdorff distance of a sequence of $r$-convex sets implies the convergence of their boundaries, as is proved in Theorem 3 in \cite{cuevas:12}, which, in turn, implies  the convergence in measure, i.e.  $d_\mu(C_{r}(S_T^{\mathrm{ON}}),S)\to 0$ \ a.s.  The choice of $r$ is important in practice, although it can be underestimated. On this respect a data-driven method is proposed in \cite{rodsav22}. Observe also that the condition $S_T^{\textrm{ON}}\subset S$ is fulfilled when $S_T^{\textrm{ON}}$ is the $r$-convex hull of $X_T^{\textrm{ON}}$ and $S$ is $r$-convex.
	\end{remark}
	
	\subsection{Estimation of the stationary distribution, level sets, and drift}
	
	In the stochastic differential equation (\ref{sde}), the drift $\nu(x)$ is given by the gradient $\nabla f$ of a function $f$, i.e.
	$\nu(x) =  -\frac{1}{2} \nabla f (x)$.
	The density of the stationary distribution (for the RBMD and the on--off process $X^{ON}_{T}$) is given by $\pi(dx) = c_0e^{ -f (x)} \mathbb{I}_{\textnormal{int}(S)}(x) dx := g(x)dx$, where $c_0$ is a normalization constant, as is stated in the following proposition.
	
	\begin{proposition} \label{sdist}
		Let $S$ satisfy $S=\overline{\textnormal{int}(S)}$, and suppose $\textnormal{int}(S)$ is  a bounded domain and  $\partial S$ is $C^2$. 
		Assume that $\nabla f$ is Lipschitz on $S$.
		Then 	
		the measure
		\begin{equation*} 
			\pi(dx)=c_0e^{-f(x)}\mathbb{I}_{\textnormal{int}(S)}dx
		\end{equation*}
		is the unique stationary measure of $\{X_t\}_{t\geq 0}$.
	\end{proposition}
	
	If $\nu=0$, the Reflected Brownian Motion without drift is obtained and the stationary distribution is the uniform distribution in $S$, \citep{burdzy:06}.  For any (smooth enough) function $\nu$, the density $g$ can also be estimated using a kernel-based estimator 
	\begin{equation}\label{kernel}
		\hat{g}_n(x)=\frac{1}{n\tau^d_n}\sum_{i=1}^n K\Big(\frac{x-X_i}{\tau_n}\Big){\ind_{C_r(X_T^{\mathrm{ON}})}(x)},
	\end{equation}
	where $K:\mathbb{R}^d\to \mathbb{R}$ is a non-negative function, as proved in \cite{ch:20}, based on a subsample $\aleph_n=\{X_1,\dots,X_n\}$ of points at which the GPS is on.
	
 It is easy to prove, following the ideas used to prove Proposition 2 in \cite{ch:20}, that the chain  $\aleph_n=\{X_{(k+1)\delta_1+k\delta_2}:k=0,\dots,n-1\}$ is geometrically ergodic. Then, from Corollary 1 together with Remark 4 in \cite{ch:20} we obtain the following remark

	\begin{remark}\label{teo:levels}
		Under the hypotheses of Theorem \ref{onoffTEO}, assume further that $g$ is Lipschitz.
		Let $\hat{g}_n$ be given by \eqref{kernel}, based on $\aleph_n$.
		Assume that $K$ is non-negative and Lipschitz, and that $\int K(t)dt=1$.
		Let $\tau_n\to 0$ and  $\gamma_n\to \infty$, such that $\gamma_n \tau_n\to 0$ and  $\log(n)/\gamma_n\to 0$.
		Then
		$$\gamma_n\sup_{x\in S} |\hat{g}_n(x)-g(x)|\to 0 \quad \mathrm{a.s.}$$
		Moreover, if $\lambda>0$ is  such that $\partial G_g(\lambda)\neq \emptyset$, where  $G_g(\lambda)=\{g>\lambda\}$,  while $g$ is $C^2$ on a neighbourhood $E$ of   $\{g=\lambda\}$ and the norm of the gradient of $g$ is strictly positive on $E$, then
		$d_H(\partial G_g(\lambda),\partial G_{\hat{g}_n}(\lambda))=o(1/\gamma_n)\quad$ a.s.
	\end{remark}

	The subsample considered in \Cref{teo:levels} can be replaced by any subsample $\aleph_n=\{X_{t_1},\dots,X_{t_n}\}$ where $t_i$ is such that the process is observed at time $t_i$, whenever $\aleph_n$ is  geometrically ergodic.
	This last condition is guaranteed for instance if there exists $\rho>0$ (independent of $n$) such that for all $i=1,\dots,n$, $t_i-t>\rho$, where $t$ is the first observed time previous to $t_i$.

	The level sets will provide significant information about the time spent in those regions, in particular the core-area will correspond to level sets with large values of $\lambda$.
	
	An estimator of the drift function can also be derived from a  plug-in method, and is given by
	$
	\hat \nu(x)= \frac{1}{2}\nabla \log(\hat g_n(x))
	$.

	\section{Some simulation results and an example with real data}
 
	Through some simulation examples and an example of real data, we compare the performance of the on--off model  with the full model (i.e. where the GPS is always on). 
 
	The simulations and computational calculations were performed using the \textbf{R} programming language \citep{R}, utilizing integration with \textbf{C++} provided by \textbf{Rcpp}. For the computation of  the $r$-convex hull , we employed the \textbf{RcppAlphahull} library \citep{rcppalphahull}, which offers faster versions of functions from the \textbf{alphahull} library \citep{alphahull}. Custom \textbf{Rcpp} code was developed for simulating RBMD and calculating distances, and it is available in the public repository at \url{https://github.com/emehache/RBM}.

	To illustrate Theorem \ref{onoffTEO}, we assume $\delta_1=\delta_2$ in subsection \ref{sec:distance_estimation}.
	On the other hand, to gain some insight about some other possible choices for $\delta_1$ and $\delta_2$ (see Remark \ref{losdelta}), we assume in the level-set estimation problem $\delta_1\neq\delta_2$, see subsection \ref{levesti}.
	This is also the case for the example with real data, in subsection \ref{realdata}.

\begin{figure}[ht]%
		\begin{center}
			\includegraphics[scale=.6]{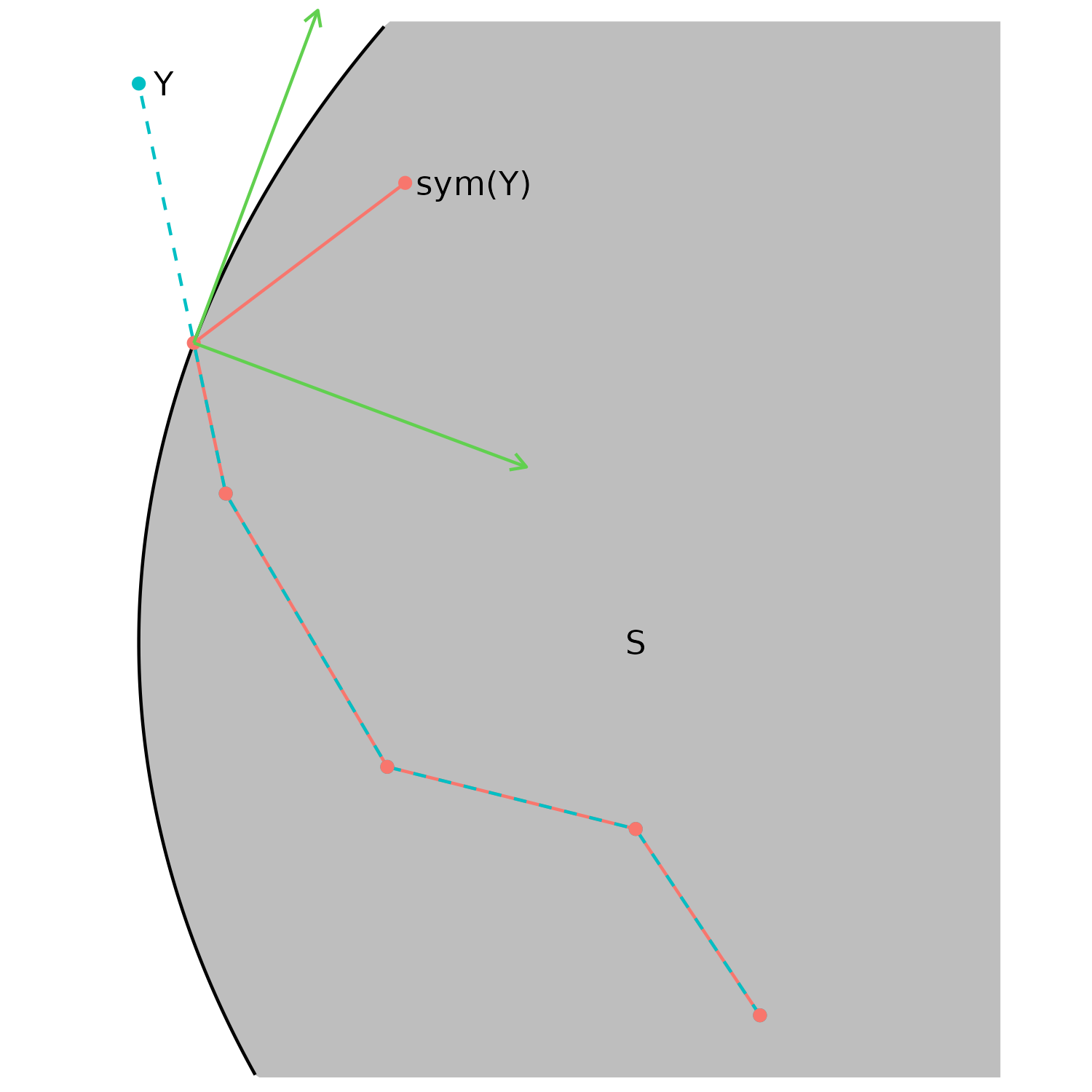} 
			\caption{An example of  the simulation scheme defined by (1), (2) and (3).}
			\label{fig:simetrizacion}
		\end{center}
	\end{figure}

	To simulate the RBMD, we follow  \cite{ch:20}: we first choose a step  $h>0$,   
	and denote by $\operatorname{sym}(z)$ the point symmetric
	to the point $z$ with respect to $\partial S$ (see Figure \ref{fig:simetrizacion}).
 
	We start with $X_0=x$ and suppose that we have obtained $X_i\in S$.
	To produce the following point, set
	\[
	Y_{i+1}=X_i+Z_i+h \nabla f(X_i),
	\] 
	where $Z_i$ is a centred Gaussian random vector,  independent of $Z_1,\dots,Z_{i-1}$, 
	with covariance matrix 
	$h \times (I_d)_{\mathbb{R}^2}$.
	
	Then:
	\begin{enumerate}
		\item If $Y_{i+1}\in S$, set $X_{i+1}=Y_{i+1}$.
		\item If $Y_{i+1}\notin S$ and $\operatorname{sym}(Y_{i+1})\in S$, set $X_{i+1}=\operatorname{sym}(Y_{i+1})$.
		\item If $Y_{i+1}\notin S$ and $\operatorname{sym}(Y_{i+1})\notin S$, set $X_{i+1}=X_i$.
	\end{enumerate}
	
In \Cref{fig:simetrizacion} we illustrate this simulation scheme.
	
	Lastly, the on--off model is obtained from $X_1,\dots,X_N$ where we only keep those $X_i$ such that $i \in \cup_{k=0}^\infty \{[k(\delta_1 + \delta_2)/h, (k+1)\delta_1/h + k\delta_2/h]\}$.
	
	We consider an RBMD in the set 
	$S=E\setminus B((4/5,0),1/2)$, where 
	$E=\{(x,y)\in \mathbb{R}^2\colon 4x^2/9+y^2\leq 1\}$, with drift function given by  
	$\nu(x,y)=-(x,y)$.
	The stationary density is
	\begin{equation}\label{case1f} 
		g(x)=c_0e^{-(x^2+y^2)}\mathbb{I}_S(x,y)\quad \text{where } c_0^{-1}=\iint_S \exp\Big[-(x^2+y^2)\Big]dxdy.
	\end{equation}
	
	\subsection{Distance estimations}
	\label{sec:distance_estimation}
	
	To show the performance of our estimator we consider 500 replications. We calculate the Hausdorff distance  between $S$ and three competitors:  the full trajectory observed on $T=p\delta_1+(p-1)\delta_2$;  the on--off model for $\delta_1=\delta_2$; and the trajectory observed on $[0,p\delta_1]$. We also compute the distance in measure between $S$ and the $0.4$-convex hull of the trajectories observed under each scenario. We choose  $h=0.0005$, $h=0.001$ as two different discretization steps, and the whole trajectories are of length $N=10^5$. 
	
	As can be seen, there is a  performance improvement when we use the on--off model (which uses half of the battery life of the full trajectory) with respect to the process observed only on $[0,p\delta_1]$. On the other hand, the loss with respect to the full trajectory is negligible. Results are shown in \Cref{tab:distancias}.

\begin{figure}[ht]%
		\begin{center}
			\includegraphics[width=.8\textwidth]{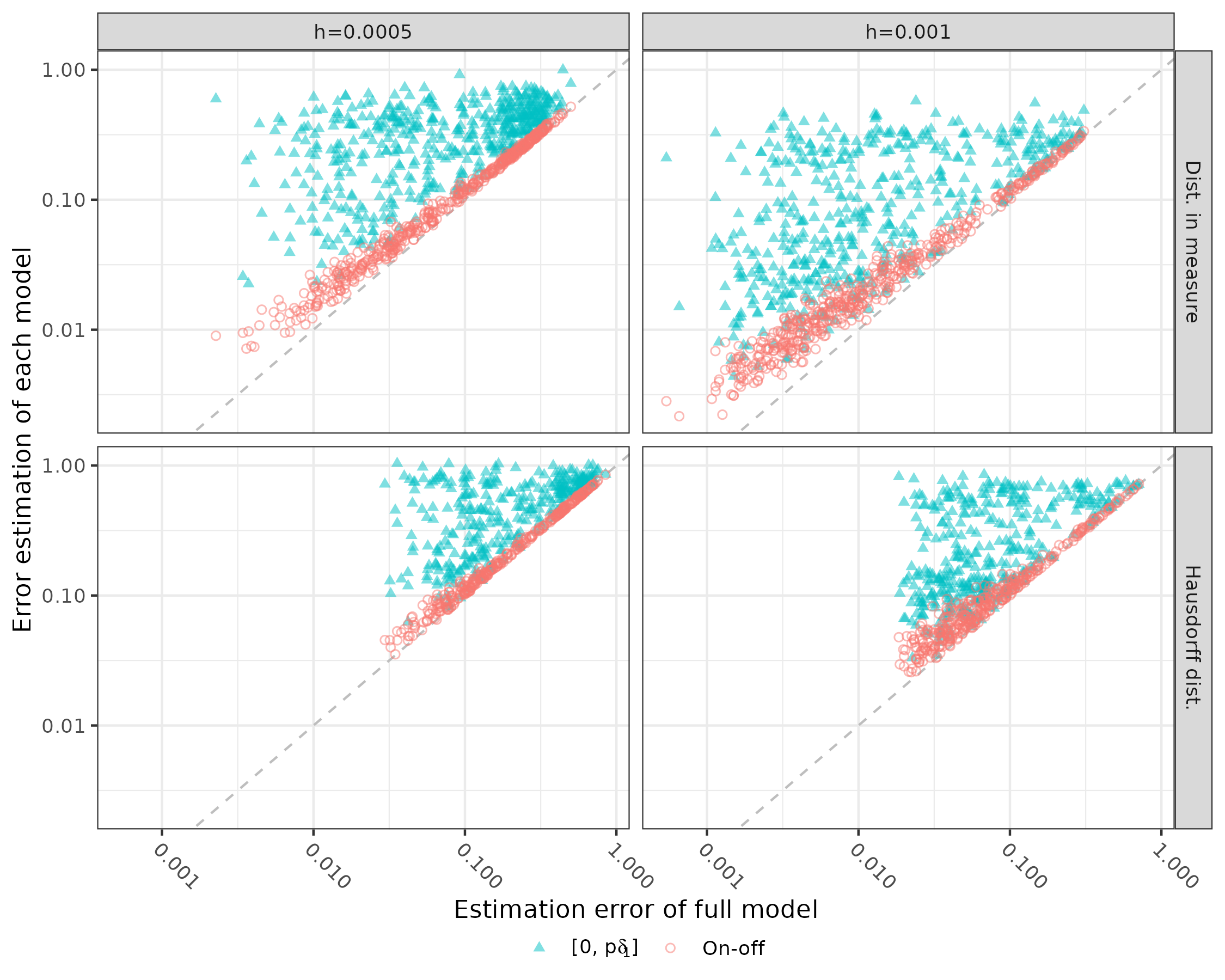} 
			\caption{				
				Points with coordinates
				$(d_H(X_T),S), d_H(X_{p\delta_1}, S))$
				and $(d_H(X_T,S),d_H(X_T^{\mathrm{ON}}, S))$ (green triangles and red circles on first row respectively) and  $(d_\mu(C_{0.4}({X_T}),S), d_\mu(C_{0.4}(X_{p\delta_1}), S))$
				and $(d_\mu(C_{0.4}({X_T}),S),d_\mu(C_{0.4}(X_T^{\mathrm{ON}}), S))$ (green triangles and red circles on second row respectively).
				In gray dashed line we plot $y=x$.
				$N=10^5$, $p=10^4$, $\delta_1/h = \delta_2/h= 5$.
			}
			\label{fig:distancias}
		\end{center}
	\end{figure}

	To illustrate the comparison, in \Cref{fig:distancias} we plot, for different replicas, points with $x$-coordinate given by the estimation error considering the full trajectory, and $y$-coordinate given by the estimation error under the different models (the on--off, and the trajectory observed on $[0, p\delta_1]$). We do this for both distances considered and two different values of the discretization step $h \in \{0.0005, 0.001\}$.

	As it can be seen, the red dots are closer to the diagonal, which shows that the estimator with the on-off trajectory is much closer (in terms of Hausdorff distance and measure) to the estimation with the full trajectory observed.

	\begin{table}
\begin{tabular}{llrrrrrr}
	\toprule
	\multicolumn{1}{c}{ } & \multicolumn{1}{c}{ } & \multicolumn{3}{c}{Hausdorff distance} & \multicolumn{3}{c}{Distance in measure} \\
	\cmidrule(l{3pt}r{3pt}){3-5} \cmidrule(l{3pt}r{3pt}){6-8}
	$h$ &   & Full & $[0,p\delta_1]$ & On--off & Full & $[0,p\delta_1]$ & On--off\\
	\midrule
	0.0005 & Mean & 0.2731 & 0.4956 & 0.2815 & 0.1266 & 0.3425 & 0.1402\\
	0.0005 & Median & 0.1828 & 0.4825 & 0.1882 & 0.0900 & 0.3440 & 0.1051\\
	0.001 & Mean & 0.1270 & 0.2932 & 0.1396 & 0.0397 & 0.1410 & 0.0485\\
	0.001 & Median & 0.0689 & 0.1864 & 0.0856 & 0.0092 & 0.0946 & 0.0172\\
	\bottomrule
\end{tabular}
		\caption{Mean and median estimation errors were computed over 500 replications for three scenarios: when considering the full trajectories, the trajectories observed under the on-off model, and the ones observed within the interval $[0, p\delta_1]$, $N=10^5$, $p=10^4$, and $\delta_1/h = \delta_2/h= 5$.}\label{tab:distancias}
	\end{table}

	\subsection{Level set estimation}\label{levesti}
	
\begin{figure}[ht]%
		\begin{center}
			\includegraphics[scale=.6]{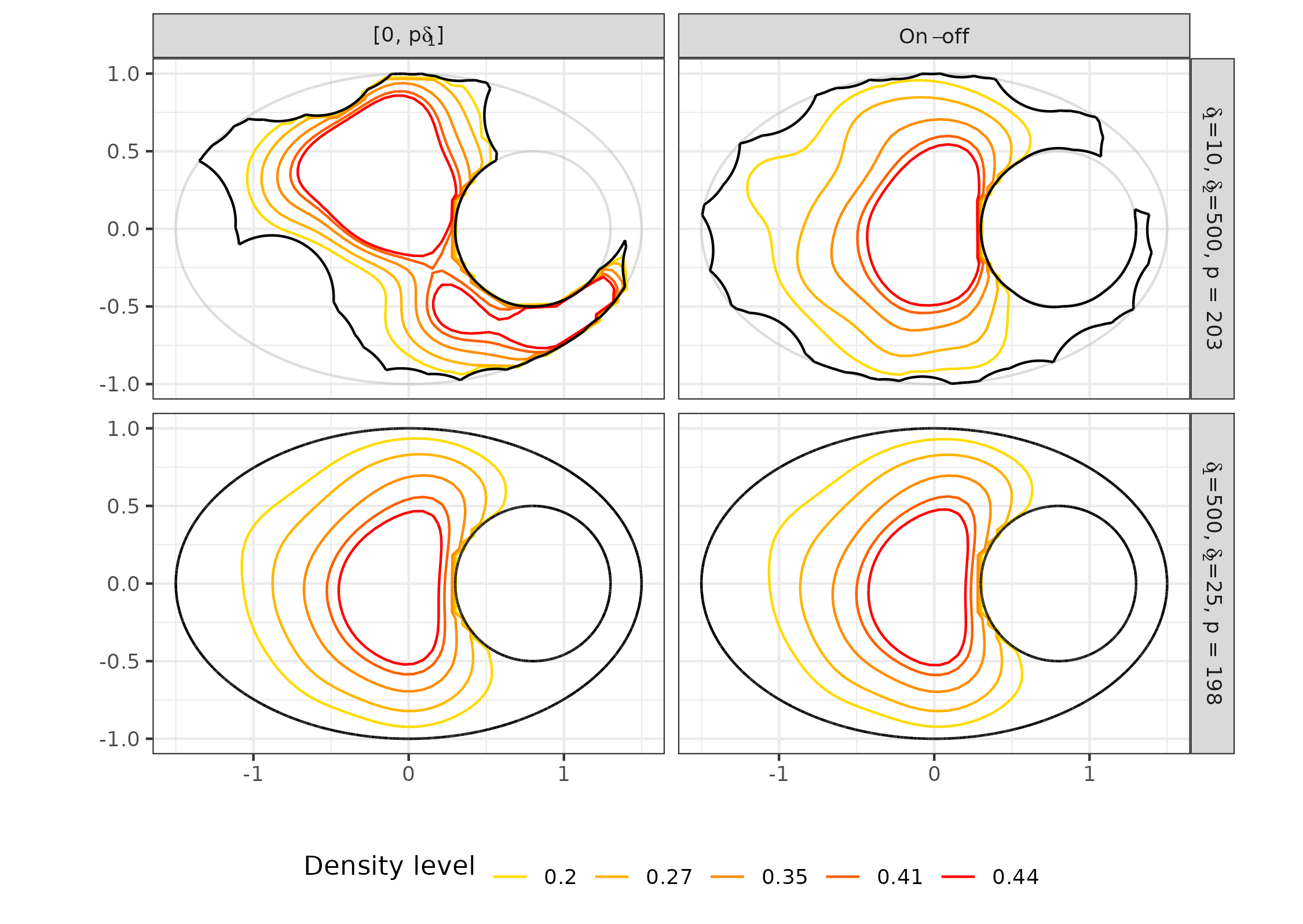} 
			\caption{Top: Contour plot of the level sets of the estimated density function when $\delta_1 = 10$, $\delta_2 = 500$ with a total number of observations of $p \delta_1 = 2030$.
				Bottom: Contour plot of the level sets of the estimated density function when $\delta_1 = 500$, $\delta_2 = 10$ with a total number of observations of $p \delta_1 = 98809$.
				Left panels correspond to the entire trajectory, right panels to the on--off model.}
			\label{fig:cont}
		\end{center}
	\end{figure}

\begin{figure}[ht]%
		\begin{center}
			\includegraphics[scale=.4]{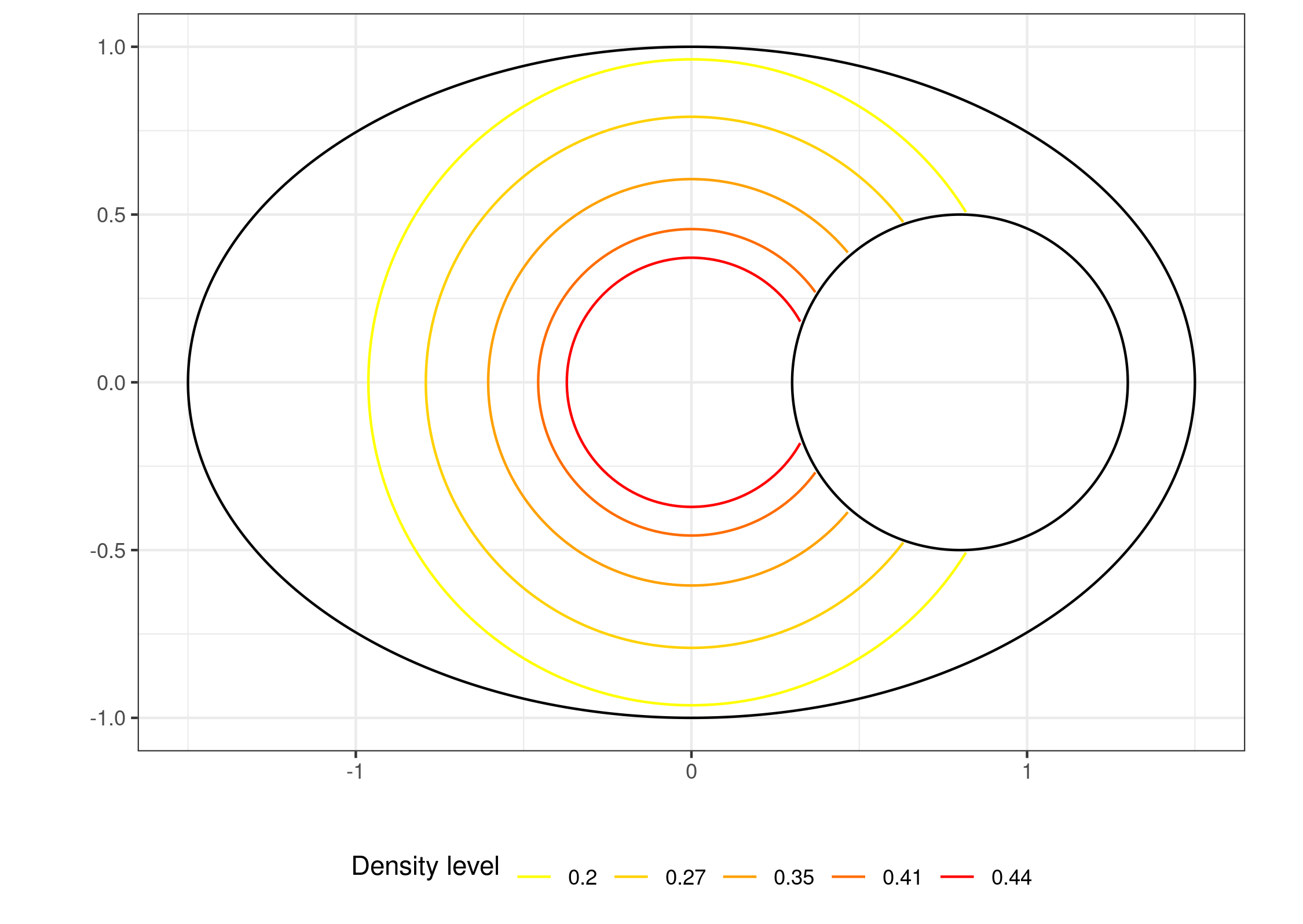} 
			\caption{Theoretical level sets for the density $g$ given by \eqref{case1f}.}
			\label{teo}
		\end{center}
	\end{figure}

	Figure \ref{fig:cont} shows the level sets of the estimator \eqref{kernel} for two different choices of $\delta_1$ and $\delta_2$, where a Gaussian kernel is used, with bandwidth $\tau=0.2$ selected by cross-validation.
	The theoretical level sets are shown in Figure \ref{teo}.
	The much better behaviour of the on--off model is clear when the number of points in the trajectory is small (2030 in the top panels), while the behaviour becomes similar when this number is large (98809 in the bottom panel).


	\subsection{An example using real data} \label{realdata}
	
\begin{figure}[h]%
		\begin{center}
			
			\makebox[\textwidth][c]{
				\includegraphics[width=1\textwidth]{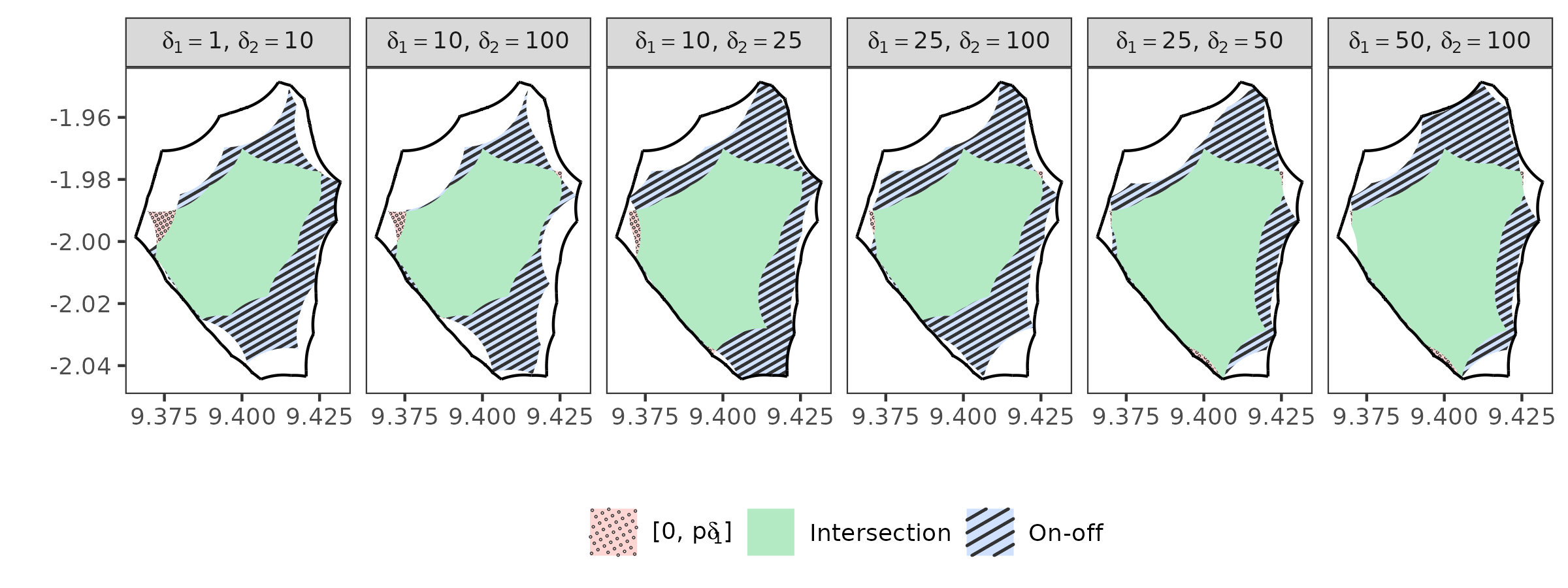}
			}%

			\caption{Each panel shows, for different values of $\delta_1$ and $\delta_2$, the $0.02$-convex hull of the trajectory observed under the on--off model compared with the $0.02$-convex hull of that observed without interruption for the same length of time. The intersection of both sets is shown in green, and the differences are shown in red and blue. In black solid line, the boundary of the $0.02$-convex hull of the whole trajectory available.}
			\label{fig5}
		\end{center}
	\end{figure}

	In this section we demonstrate the performance of the on--off model using an example of real data.
	We consider a data set consisting of 1633 recorded positions of elephants in Loango National Park in western Gabon. The data that support the findings of this study are openly available in Movebank at \url{https://www.movebank.org/cms/webapp?gwt_fragment=page=studies,path=study1818825}, reference number 1818825.
	This dataset was also analysed in \cite{ch:20}.
	We first estimate the $r$-convex hull of this full trajectory.
	Later, we imagine that this full trajectory is not available at all, and we only have a subset of size $p\delta_1$ of the recorded locations.
	One approach is to observe the first $p\delta_1$ steps, and the other approach is to use our on--off strategy.
	Figure \ref{fig5} shows, as a solid black line, the boundary of the $0.02$-convex hull of the full trajectory, and the $0.02$-convex hull under the two approaches. The intersection of the two is shown in green, and the differences between one and the other are shown in red and blue.  We assume that the $0.02$-convex hull is a good approximation of the true home-range. Then, Figure \ref{fig5} shows that the $0.02$-convex hull of the on-off trajectory performs much better than the $0.02$-convex hull of the trajectory observed between $0$ and $p\delta_1$. 
	
\begin{figure}[h]%
		\begin{center}
			\makebox[\textwidth][c]{
				\includegraphics[width=1\textwidth]{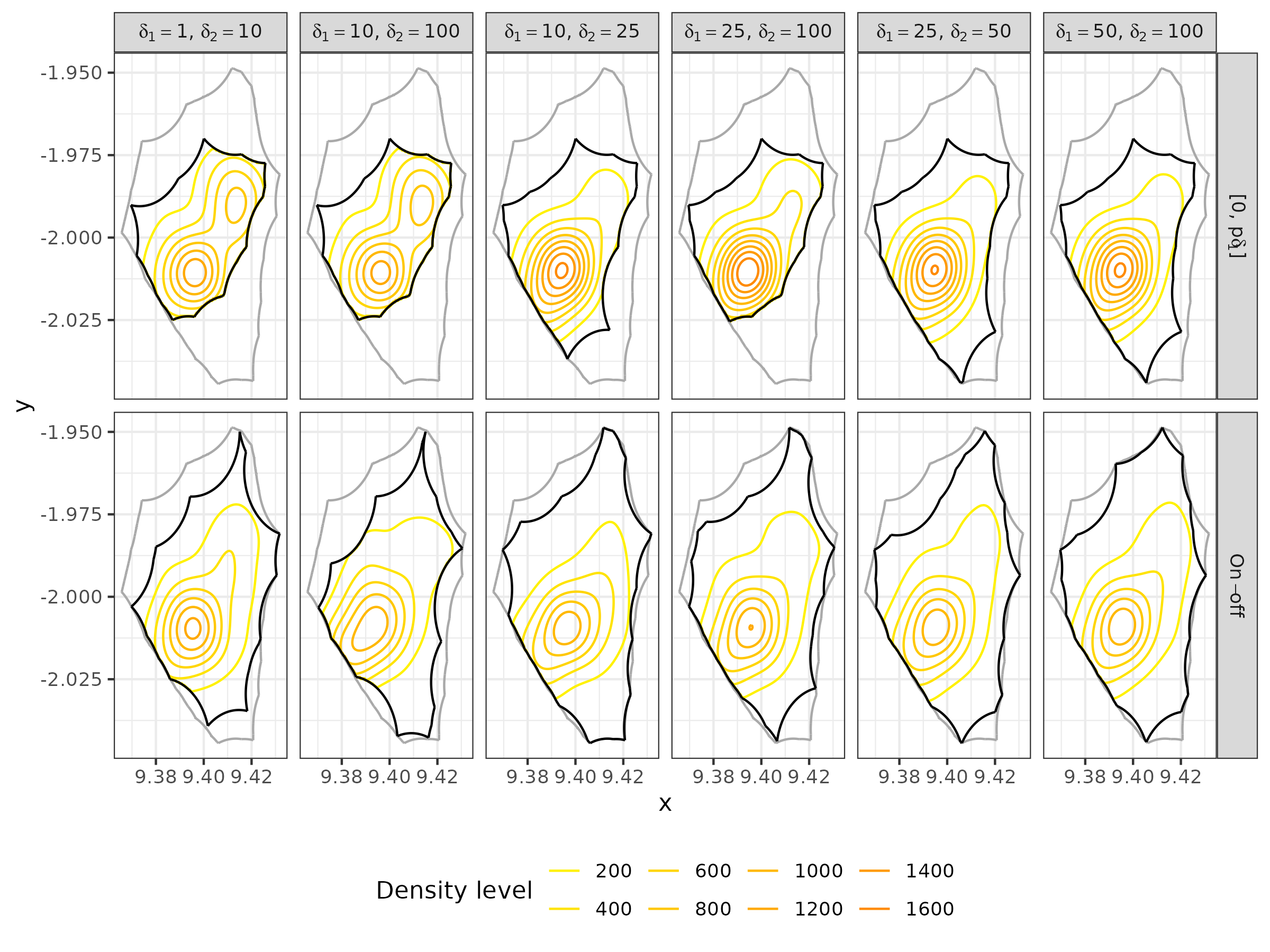}
			}%
 
			\caption{Contour plot of the level sets of the estimated density function for different choices of $\delta_1$ and $\delta_2$. The gray line represents the $0.02$-convex hull of the entire available trajectory. The black line represents the counterpart under each model. Top: cases where the trajectory is continuously observed on $[0, p\delta_1]$. Bottom: cases where the trajectory is observed using the on-off model. We used a Gaussian kernel with a bandwidth parameter of $\tau = 0.028$.}
			\label{fig6}
		\end{center}
	\end{figure}

	The stationary distribution is also estimated. Figure \ref{fig6} presents the level set estimation for the probability density function under different choices of the parameters $\delta_1, \delta_2$.

	\section{Concluding Remarks}
	
	We obtain almost sure consistency results for the estimation under the Hausdorff distance and the distance in measure for the on--off model. We also obtain the rate of convergence  for the Hausdorff distance, which turns out to be almost optimal (i.e. equal, up to a logarithmic term, to the case in which the data is an iid sample). The same rate of convergence for the Hausdorff distance is obtained if the trajectory is observed by a GPS that is always on.

	The stationary distribution is estimated using a kernel type estimator as proposed in \cite{ch:20}.
	
	From the uniform convergence of the estimated stationary distribution, we derive estimators of the level sets, which can determine the core-area of the animals' home-range. 
	
	An estimator of the drift function can be derived from an estimator of the stationary distribution by a simple plug-in rule.
	
	An optimal choice of the parameters $\delta_1, \delta_2$ remains an open problem.
	
	It could be interesting to generalize the on--off model in future work, allowing $\delta_1$ and $\delta_2$ to be random variables, so the sampling scheme would also be a random variable. Such a new setup could include the case when an accelerometer governs the power status of the GPS, even also some missing data due to connection problems could be admitted.

	\section*{Acknowledgements} 		
	We thanks Dr. Stephen Blake, of the Max Planck Institute for Ornithology, for facilitating access to the data set that was used in this paper.
	The data that support the findings of this study are openly available in Movebank at \url{https://www.movebank.org/cms/webapp?gwt_fragment=page=studies,path=study1818825}, reference number 1818825.
	
 We thanks the editor and three referee’s for their constructive comments which improves significantly the present version of the manuscript.
	
	This work was supported by grants FCE120191156054, ANII and POSNAC20191157608, ANII.

	\section{Appendix}

	\subsection{Proof of Theorem  \ref{onoffTEO} a)}
	
	Assume first that $\delta_1$ and $\delta_2$ are positive and fixed. Let us fix $\eps>0$.   We will define a grid depending on $\eps$, on $[0,T]\cap \mathcal{I}$.  Let us define $t_\eps:=-(1/\alpha)\log(\omega_dc (\eps/2)^d/(2\beta ))$.  
	Let  us define $t_i:= i\times s(\delta_1+\delta_2)$ for $i\in \mathbb{N}$, where $s$ is the smallest integer that guarantees $s(\delta_1+\delta_2)>t_\eps$.  Then, for all $n$, there exists $T$ large enough, such that there exist
	$t_1,\dots,t_n \in \mathcal{I}\cap [0,T]$ and $t_{i+1}-t_i>t_\eps$ for all $i=1,\dots,n-1$. We will assume that the process is ON in $[T-\delta_1,T]$.  
	
	Denote the $\eps-$inner parallel set of $S$ by $S^{(\eps)}=\{x \in S: B(x,\eps) \subset S\}.$
	
Put $I_n:=\{1,\dots,n\}$, observe that $a_{t_i}=1$ for all $i\in I_n$. Then \begin{align*}
		\mathbb{P}\{d_H(S_T^{\mathrm{ON}},S)>\eps \} \leq & \ \mathbb{P}\{ \exists x \in S^{(\eps)}: \forall t \in \mathcal{I},t<T : X_{t} \not\in B(x,\eps)\} \\
		\leq &\  \mathbb{P}\{ \exists x \in S^{(\eps)}: \forall i \in I_n: X_{t_i} \not\in B(x,\eps)\}.
	\end{align*}	
	Let $x_1,\dots,x_N \in S^{(\eps)}$ be such that $S^{(\eps)} \subset B(x_1,\eps/2) \cup \dots \cup B(x_N,\eps/2)$, and suppose that $N$ is the smallest positive integer such that such a covering of $S^{(\eps)}$ is possible. Then $N=N(\eps /2)$ is called the $\eps /2$-covering number of $S^{(\eps)}$. Since $B(x_i,\eps/4)\subset S$ for all $i=1,\dots,N$ and $B(x_i,\eps/4)\cap B(x_j,\eps/4)$ for all $i\neq j$, $N \leq \mu(S)/\mu(B(0,\eps/4))=(\eps/4)^{-d}\mu(S)/\omega_d $.
	
	If for some $x \in S^{(\eps)}$ we have $X_{t_i} \not\in B(x,\eps)$ for all $i \in I_n$, then there exists a $j \in \{1,\dots,N\}$ such that $X_{t_i} \not\in B(x_j,\eps/2)$ for all $i=1,\dots,n.$ Thus, continuing the chain of inequalities above, 
	\begin{align*}
		\mathbb{P}\{d_H(S_T^{\mathrm{ON}},S)>\eps \} \leq & \ 
		\mathbb{P}\{ \exists j \in \{1,\dots,N\}: \forall i \in I_n: X_{t_i} \not\in B(x_j,\eps/2)\} \\
		\leq & \  N \sup_{x\in S^{(\eps)}} \mathbb{P}\{\forall i \in I_n: X_{t_i} \not\in B(x,\eps/2)\}.
	\end{align*}
	Next we estimate the probability on the right-hand side. Recall that the process is ON in $[T-\delta_1,T]$, and so $n\in I_n$. For all $x \in S^{(\eps)}$,
	\begin{align*}
		\mathbb{P}\{\forall i \in I_n: X_{t_i} \not\in B(x,\eps/2)\} \\
		= & \ \mathbb{P}\{X_{t_n} \not\in B(x,\eps/2)|\forall i \in I_{n-1}: X_{t_i} \not\in B(x,\eps/2) \} \\ 
		&\ \times \mathbb{P}\{\forall i \in I_{n-1}: X_{t_i} \not\in B(x,\eps/2)\} \\
		= & \ \mathbb{P}\{X_{t_n} \not\in B(x,\eps/2)|X_{t_{n-1}} \not\in B(x,\eps/2) \} \\
		& \  \times \mathbb{P}\{\forall i \in I_{n-1}: X_{t_i} \not\in B(x,\eps/2)\} \\
		& \ \text{(since $X_t$ is a Markov process)}
	\end{align*}
	Let us iterate this process
	
	$$	\mathbb{P}\{\forall i \in I_n: X_{t_i} \not\in B(x,\eps/2)\} = \prod_{i=0}^{n-1} \mathbb{P}\{X_{t_{n-i}} \not\in B(x,\eps/2)|X_{t_{n-i-1}} \not\in B(x,\eps/2) \}.$$
	
	Now, by Proposition \ref{properg}, 
	\begin{align}\label{upperbound}
		\mathbb{P}\{X_{t_{n-i}} \not\in B(x,\eps/2)|X_{t_{n-i-1}} \not\in B(x,\eps/2) \} &\nonumber\\
		&\hspace{-3cm}  =1- \mathbb{P}\{X_{t_{n-i}} \in B(x,\eps/2)|X_{t_{n-i-1}} \not\in B(x,\eps/2) \} \nonumber\\
		&\hspace{-3cm}\leq   1-\pi(B(x,\eps/2))+\beta \exp\{-\alpha(t_{n-i}-t_{n-i-1})\} \nonumber\\
		&\hspace{-3cm} =  1 - c\omega_d \eps^d/2^{d+1}.
	\end{align}
	$$\prod_{i=0}^{n-1} \mathbb{P}\Bigg\{X_{t_{n-i}} \not\in B(x,\eps/2)|X_{t_{n-i-1}} \not\in B(x,\eps/2)\Bigg \}\leq \Bigg(1-\frac{c \omega_d \eps^d}{2^{d+1}}\Bigg)^{n}\leq \exp\Bigg(\frac{-nc \omega_d \eps^d}{2^{d+1}}\Bigg)$$
	Then
	$$ 	\mathbb{P}\{d_H(S_T^{\mathrm{ON}},S)>\eps \} \leq \frac{(\eps/4)^{-d}\mu(S)}{\omega_d}   \exp\Bigg(\frac{-nc \omega_d \eps^d}{2^{d+1}}\Bigg).$$ 
	From the Borel--Cantelli Lemma, there follows a).

	\subsection{Proof of Theorem  \ref{onoffTEO} b)}
	
	The proof follows the same ideas used to prove part a).   Fix $\eps>0$.   We will define a grid depending on $\eps$ and $T$, on $[0,T]\cap \mathcal{I}$.  Let  us define, for $i\in \mathbb{N}$, $t_i= i\times (\delta_1+\delta_2)$. 
	We will assume that the process is ON in $[T-\delta_1,T]$.  Proceeding as in \eqref{upperbound}, 
	
	\begin{align*} 
		\mathbb{P}\{X_{t_{n-i}} \not\in B(x,\eps/(2\kappa_T))|X_{t_{n-i-1}} \not\in B(x,\eps/(2\kappa_T)) \}  &\\
		&\hspace{-4cm} = 1- \mathbb{P}\{X_{t_{n-i}} \in B(x,\eps/(2\kappa_T))|X_{t_{n-i-1}} \not\in B(x,\eps/(2\kappa_T)) \} \\
		&\hspace{-4cm}\leq   1-\pi(B(x,\eps/(2\kappa_T))+\beta \exp\{-\alpha(\delta_1+\delta_2)\} \\
		&\hspace{-4cm}\leq   1-c\eps^d/(2\kappa_T)^d+\beta \exp\{-\alpha(\delta_1+\delta_2)\}.\\
	\end{align*}
	Since $\eta_T\to \infty$, for all $\eps$ there exists $T$ large enough such that 
	$$\beta \exp\{-\alpha(\delta_1+\delta_2)\} \leq c\eps^d/(2^{d+1}\kappa_T^d).$$
	Since $n=\lfloor T/(\delta_1+\delta_2)\rfloor$, we get
	\begin{multline*}
		\prod_{i=0}^{n-1} \mathbb{P}\Bigg\{X_{t_{n-i}} \not\in B(x,\eps/(2\kappa_T))|X_{t_{n-i-1}} \not\in B(x,\eps/(2\kappa_T)) \Bigg\}\leq \\
		\Bigg(1-\frac{c \omega_d \eps^d}{2^{d+1}\kappa_T^d}\Bigg)^{T/[\delta_1+\delta_2]}\leq \exp\Bigg(-\frac{Tc \omega_d \eps^d}{2^{d+1}\kappa_T^d}\Bigg).
	\end{multline*}
 
	Then, for all $\eps$, for $T$ large enough,
	\begin{align*}
		\mathbb{P}\{d_H(S_T^{\mathrm{ON}},S)>\eps \}\leq & (\eps/4)^{-d}\mu(S)/\omega_d\exp\Bigg(-\frac{Tc \omega_d \eps^d}{2^{d+1}\kappa_T^d}\Bigg)\\
		=&(\eps/4)^{-d}\mu(S)/\omega_d\exp\Bigg(-\frac{\log^2(T)c \omega_d \eps^d}{2^{d+1}}\Bigg).
	\end{align*}
	Lastly, part b) follows from the Borel--Cantelli Lemma.

	\subsection{Proof of Proposition \ref{sdist}}

In order to obtain the stationary distribution, we will use the following lemma, whose proof is accomplished by reasoning as in the proof of lemma 2.1 $(i)$ in \cite{hw:87}. First let us introduce some notation.  We denote by $C_c^2(S)$ the set of twice continuously differentiable functions with compact support in some domain containing $S$. We write $\mathcal{L}$ for the infinitesimal generator of the process, i.e. $\mathcal{L}(h)(x)=\lim_{t\downarrow 0}(1/t)(\mathbb{E}_x(h(X_t))-h(x))$. It can be proved that $\gen h=\frac{1}{2}\Delta h -\frac{1}{2}\langle \nabla f, \nabla h \rangle$, for $h\in C_c^2(S)$, \citep{ch:20}.

\begin{lemma} \label{lemstdens} Let $S$ satisfy $S=\overline{\textnormal{int}(S)}$, and suppose $\textnormal{int}(S)$ is  a bounded domain and  $\partial S$ is $C^2$. Suppose that $p\colon S\rightarrow \mathbb{R}$ is $C^2$, positive on $int(S)$, and that $\int_{\textnormal{int}(S)}p(x)dx=1$. Then $p$ is the density of the (unique) invariant distribution for \eqref{sde} if and only if
	\begin{equation*}
		\int_{\textnormal{int}(S)} p(x) \gen h(x)dx=0\quad \text{ for all } h\in C_c^2(S) \text{ satisfying }\langle \nabla h(x),\normal(x)\rangle =0 \text{ on }\partial S,
	\end{equation*}
	where $\normal(x)$ denotes the inner normal vector at $x\in \partial S$.
\end{lemma}

\textit{Proof of Proposition \ref{sdist}}\\
By Lemma \ref{lemstdens}, the measure $\pi$ is the stationary distribution if and only if for all $h\in C_c^2(S)$ with $\langle \normal(x), \nabla h(x)\rangle=0$, 
for all $x\in \partial S$, one has that $0=\int_{\textnormal{int}(S)}ce^{-f(x)} \gen h(x)dx$. But this is a direct consequence of Green's first identity:
\begin{align*}
	-\int_{\textnormal{int}(S)}e^{-f(x)}\Delta h(x)=&\int_{\partial S}e^{-f(x)}\langle \nabla h(x), \normal(x)\rangle d\sigma(x)+\int_{\textnormal{int}(S)}e^{-f(x)}\langle \nabla f(x), \nabla h(x)\rangle dx\\
	=&\int_{\textnormal{int}(S)} e^{-f(x)}\langle \nabla f(x), \nabla h(x)\rangle dx,
\end{align*}
with $\sigma$ being the surface measure on $\partial S$.

\end{document}